\documentclass[a4paper,11pt,makeidx]{amsart}
\usepackage{amscd}
\usepackage{xypic}  
\usepackage{amssymb}
\usepackage{amsthm}
\usepackage{epsf}
\makeindex

\theoremstyle{plain}

\def\mult{\operatorname{mult}}

\def\c1{\operatorname{c_1}}
\def\c2{\operatorname{c_2}}

\def\ZZ{{\mathbb Z}}
\def\NN{{\mathbb N}}

\def\*{\otimes}

\def\+{\oplus}                   
\def\*{\otimes}                  

\def\Pic{\operatorname{Pic}}

\begin{document}

\title{A note on Seshadri constants on general $K3$ surfaces}

\author{Andreas Leopold Knutsen}

\address{\hskip -.43cm Andreas Leopold Knutsen, Dipartimento di Matematica, 
Universit\`a di Roma Tre, Largo San
Leonardo Murialdo 1, 00146,
Roma, Italy. }

\email{knutsen@mat.uniroma3.it}

\thanks{Research supported by a Marie Curie Intra-European Fellowship within the 6th European
Community Framework Programme}


\begin{abstract}
We prove a lower bound on the Seshadri constant $\varepsilon (L)$ on a $K3$ surface $S$ with
$\Pic S \simeq \ZZ[L]$. In particular, we obtain that $\varepsilon (L)=\alpha$ if $L^2=\alpha^2$
for an integer $\alpha$.
\end{abstract}

\maketitle

\section{Introduction and results}

Let $X$ be a smooth projective variety and $L$ be an ample line bundle
on $X$. Then the real number 
\[ \varepsilon (L,x) := \inf_{C \ni x} \frac{L.C}{\mult_x C}, \]
 introduced by Demailly \cite{dem}, is the {\it Seshadri constant of $L$ at $x \in X$} (where the infimum
is taken over all irreducible curves on $X$ passing through $x$). The {\it (global) Seshadri 
constant of $L$} is defined as
\[ \varepsilon (L) := \inf_{x \in X} \varepsilon (L,x).\]

We refer to \cite[pp.~270--303]{laz} for more background, properties and results on these constants.

The subtle point about Seshadri constants is that their exact values are
known only in a few cases and 
even on surfaces it is difficult to control them. 

It is known that the global Seshadri constant on a surface satisfies 
$\varepsilon(L) \leq \sqrt{L^2}$, cf. e.g. \cite[Rem. 1]{st}, and that $\varepsilon(L)$ is rational if
$\varepsilon(L) < \sqrt{L^2}$, cf. \cite[Lemma 3.1]{sz} or \cite[Cor. 2]{og}. 
(It is not known whether Seshadri constants are always rational, but no examples 
are known where they are irrational.)

In the case of $K3$ surfaces, Seshadri constants have only been computed for the hyperplane bundle of quartic
surfaces \cite{bau2} and in the particular case of non-globally generated ample line bundles 
\cite[Prop. 3.1]{bds}.

In this note we prove the following result:

\vspace{0,3cm}

\noindent {\bf Theorem}
{\it Let $S$ be a smooth, projective $K3$ surface with $\Pic S \simeq \ZZ[L]$.
Then either
\[ \varepsilon (L) \geq \lfloor \sqrt{L^2} \rfloor, \]
or 
\begin{equation} \label{eq:exc}
 (L^2,\varepsilon(L)) \in \Big\{ (\alpha^2+\alpha-2, \alpha-\frac{2}{\alpha+1}), 
                                   (\alpha^2+ \frac{1}{2}\alpha- \frac{1}{2},\alpha-\frac{1}{2\alpha+1}) \Big\}
\end{equation}
for some $\alpha \in \NN$. (Note that in fact $\alpha= \lfloor \sqrt{L^2} \rfloor$.)}

\vspace{0,3cm}

\noindent {\bf Remark} 
In the two exceptional cases (\ref{eq:exc}) of the theorem, the proof below shows that there 
has to exist a point $x \in S$ and an irreducible rational curve $C \in |L|$ (resp. $C \in |2L|$)
such that $C$ has an ordinary singular point of multiplicity $\alpha+1$ (resp. $2\alpha+1$) at 
$x$ and is smooth outside $x$, and $\varepsilon (L)= L.C/\mult_x C$.

By a well-known result of Chen \cite{ch}, rational curves in the primitive class of a  
{\it general} $K3$ surface in the moduli space are nodal. {\it Hence
the first exceptional case in (\ref{eq:exc}) cannot occur on a {\it general} $K3$ surface in the moduli space} (as $\alpha \geq 2$). If $\alpha=2$, so that $L^2=4$, this special case is case (b) in \cite[Theorem]{bau2}.

As one also expects that rational curves in any multiple of the primitive class on a {\it general} $K3$ surface 
are always nodal (cf. \cite[Conj. 1.2]{ch2}),
we expect that also the second exceptional case in (\ref{eq:exc}) cannot occur on a {\it general} $K3$ surface.

\vspace{0.3cm}

Since $\varepsilon(L) \leq \sqrt{L^2}$,
an immediate corollary of the theorem is the following:

\vspace{0,3cm}

\noindent {\bf Corollary}
{\it Let $S$ be a smooth, projective $K3$ surface such that $\Pic S \simeq \ZZ[L]$ with $L^2=\alpha^2$ for an integer $\alpha \geq 4$.

Then $\varepsilon (L)=\alpha$.}

\section{Proof of the theorem}

The reader will recognize the similarity of the proof of the theorem with the proofs of
\cite[Thm. 4.1]{bau} and  \cite[Prop. 1]{st}.

Set $\alpha:=\lfloor \sqrt{L^2} \rfloor$ and assume that $\varepsilon (L) <\alpha$. Then it is well-known 
(see e.g. \cite[Cor. 2]{og}) that
there is an irreducible curve $C \subset S$ and a point $x \in C$ such that
\begin{equation} \label{eq:1}
C.L < \alpha \mult_x C.
\end{equation}
Set $m:=\mult_x C$. Since a point of multiplicity $m$ causes the geometric genus of an irreducible curve to drop at least by $m \choose 2$ with respect 
to the arithmetic genus, we must have
\begin{equation} \label{eq:gen}
p_a(C) = \frac{1}{2}C^2+1 \geq {m \choose 2}=\frac{1}{2}m(m-1), 
\end{equation}
so that
\begin{equation} \label{eq:2}
m(m-1)-2 \leq C^2.
\end{equation}
We have that $C \in |nL|$ for some $n \in \NN$. From (\ref{eq:1}) we obtain
$nL^2 < m\alpha$, so that, by assumption, $n\alpha^2 < m\alpha$, whence $n\alpha  < m$. As $n\alpha \in \ZZ$ we must have
\begin{equation} \label{eq:3}
n\alpha \leq m-1.
\end{equation}
Combining (\ref{eq:1}), (\ref{eq:2}) and (\ref{eq:3}), we obtain
\[ m(m-1)-2 \leq C^2 = nC.L < n\alpha m \leq m(m-1), \]
giving the only possibilities $C^2=n^2L^2=m(m-1)-2$ and $n\alpha=m-1$.
It follows from (\ref{eq:gen}) that $C$ is a rational curve with a single singular point $x$ that is an
ordinary singularity of multiplicity $m \geq 2$. 

As
\begin{equation}  \label{eq:4} 
C.L =nL^2 = \frac{m(m-1)-2}{n} =m\alpha -\frac{2}{n}
\end{equation}  
and $m\alpha \in \ZZ$, we must have
$\frac{2}{n} \in \ZZ$, so that $n=1$ or $2$.

If $n=1$, then $m=\alpha+1$, so that $L^2=C^2=m(m-1)-2=\alpha(\alpha+1)-2$ and 
$\varepsilon (L)=C.L/m=  \alpha-\frac{2}{\alpha+1}$ from (\ref{eq:4}).

If $n=2$, then $m=2\alpha+1$, so that 
$L^2=\frac{1}{4}C^2=\frac{1}{4}((2\alpha+1)2\alpha-2)$
and $\varepsilon (L) = \alpha-\frac{1}{2\alpha+1}$ from (\ref{eq:4}).

This concludes the proof of the theorem.

\end{document}